\documentclass[11pt]{article}
\usepackage{graphicx}
\graphicspath{{images/}{../images/}}
\usepackage[T1]{fontenc}
\usepackage[utf8]{inputenc}
\usepackage{amssymb}
\usepackage{amsfonts}
\usepackage{dsfont}
\usepackage{mathtools}
\usepackage{amsthm}
\usepackage{amsmath}
\usepackage{textcomp}
\usepackage{xcolor}
\usepackage{color,soul}
\usepackage{eurosym}
\usepackage{chngcntr}
\usepackage[numbers]{natbib}
\usepackage{cleveref}
\usepackage[english=british]{csquotes}
\usepackage{enumitem}
\counterwithout{table}{section}
\counterwithout{table}{subsection}

\usepackage{geometry}
\newtheorem{lemme}{Lemma}
\newtheorem{thm}{Theorem}
\newtheorem{prop}{Proposition}

\newtheorem{rem}{Remark}

\begin{document}

\title{Optimal incentives in a limit order book: a SPDE control approach\footnote{This work benefits from the financial support of the Chaire Deep Finance and Statistics. The authors would like to thank René Aïd (Université Paris-Dauphine) for the insightful discussions they had on the subject.}}

\author{Bastien \textsc{Baldacci}\footnote{Quantitative Advisory Solutions, bastien.baldacci.qas@protonmail.com.} \and Philippe \textsc{Bergault}\footnote{\'Ecole Polytechnique, CMAP, 91128, Palaiseau, France, philippe.bergault@polytechnique.edu. \textit{(Corresponding author.)}} }
\date{}
\maketitle

\begin{abstract}
With the fragmentation of electronic markets, exchanges are now competing in order to attract trading activity on their platform. Consequently, they developed several regulatory tools to control liquidity provision / consumption on their liquidity pool. In this paper, we study the problem of an exchange using incentives in order to increase market liquidity. We model the limit order book as the solution of a stochastic partial differential equation (SPDE) as in \cite{cont2021stochastic}. The incentives proposed to the market participants are functions of the time and the distance of their limit order to the mid-price. We formulate the control problem of the exchange who wishes to modify the shape of the order book by increasing the volume at specific limits. Due to the particular nature of the SPDE control problem, we are able to characterize the solution with a classic Feynman-Kac representation theorem. Moreover, when studying the asymptotic behavior of the solution, a specific penalty function enables the exchange to obtain closed-form incentives at each limit of the order book. We study numerically the form of the incentives and their impact on the shape of the order book, and analyze the sensitivity of the incentives to the market parameters.
\end{abstract}

\vspace{9mm}

\setlength\parindent{0pt}


\vspace{6mm}

\section{Introduction}

At the end of the twentieth century, the electronification of financial markets started first with stock exchanges, before spreading to every asset class. On most order-driven markets, computers now handle almost all the market activity, organized around all-to-all limit order books. The recent years have also seen the fragmentation of electronic markets, meaning that the same asset can be traded in several platforms, thus increasing the competition between exchanges to attract trading activity on their liquidity pools, see for example \cite{laruelle2018market}. The usual mechanism put in place by exchanges is the make-take fees system: a fee rebate is associated with executed limit orders while transaction costs are applied to market orders.\footnote{An other way for the exchange to attract liquidity on its platform is to set an appropriate tick size, see \cite{baldacci2020bid, dayri2015large} for instance.} This way, exchanges charge in an asymmetric way liquidity provision
and liquidity consumption. One unintended effect of these incentives mechanism has been the rise of high-frequency trading activity, see among others \cite{arnuk2010concept}, which is now responsible for the majority of liquidity provision on electronic markets. Indeed, traditional market making companies are now competing with all market participants who can post liquidity-providing orders. One of the main concerns regarding this activity is the quality of liquidity provided in times of stress, when the high-frequency traders tend to leave the market, see \cite{bellia2017high, megarbane2017behavior,menkveld2013high,mirrlees2013strategies}. The influence of the make-take fees mechanism has been studied mostly from an empirical and data-driven viewpoint, see \cite{angel2011equity,brolley2013informed,colliard2012trading,harris1982asymmetric} showing the dependence of depths, volumes, or price impact of the order book on the make-take fees structure.\\

The recent work of \cite{el2018optimal} proposes a make-take fees design, combining the Avellaneda-Stoikov \cite{avellaneda2008high} model and a Principal-Agent framework with drift control,\footnote{There is an extensive literature on this topic, see the pioneering work of \cite{sannikov2008continuous} and recent developments such as  \cite{elie2019contracting,elie2019tale}.} to increase market liquidity on a single asset. In their formulation, the exchange offers a contract, resulting from a Stackelberg equilibrium between the market maker and the exchange, to the market maker in order to decrease the bid-ask spread offered on the asset. The optimal contract proposed is a remuneration indexed on the number of transactions realized on the asset as well as its dynamics. This contract is obtained in quasi closed-form, depending on the market maker’s inventory trajectory and market volatility. It is shown numerically that the contract effectively reduces the spreads and subsequently increases liquidity provision. This work has been extended in various ways, taking into account several market participants, the specificities of options market and enabling trading on dark liquidity pools, see \cite{baldacci2019market,baldacci2019optimal,baldacci2019design}. However, the very design of make-take fees through a Principal-Agent has three main practical drawbacks. First, the form of the incentives is a Markovian function of the inventory of the market maker. While it is reasonable for the exchange to track the inventory of the agent on the platform, it neglects the fact that she may have submitted buy or sell orders on other platforms where the asset is listed. It induces a possible arbitrage opportunity for the market maker, who can take advantage of the fact that the exchange can only see the portion of her inventory on its platform to earn a higher incentive. Second, the market making model used in these articles is the one of Avellaneda and Stoikov \cite{avellaneda2008high}, which is far more suited to Over-The-Counter (OTC) markets than to limit order books.\footnote{See \cite{bergault2021size} or \cite{gueant2016financial} for more details.} In particular, there is no notion of discrete price limits. Third, the use of agent-based models for incentives design introduces the notion of risk-aversion for the participants. While the risk-aversion parameter of the exchange can be chosen in order to produce a certain type of incentives, the one of the agent is completely unknown to the exchange. As the values of the incentives are highly dependent on this parameter, its under or overestimation can lead to very different contracts.\\ 

In this article, we propose an incentives design with a complete different paradigm. We do not rely on an agent-based model, in the sense that the incentives do not emerge from the interactions between market participants and the exchange. We consider a single asset whose order book dynamics is governed by a system of stochastic partial differential equations (SPDE) on the bid and ask sides, described in the recent work of Cont and Müller \cite{cont2021stochastic}. Two great advantages of this model are its accuracy with respect to market order-flow dynamics and its tractability via a parametrization using low-dimensional diffusion processes. We consider that the exchange is looking for the best incentives policy to offer to liquidity providers in order to increase the rate of buy or sell order submissions, so that the shape of the order book is modified. Ideally, the exchange wishes to have the maximum possible volume on the bid and ask sides at the first limits. Its willingness to increase the market liquidity decreases with the depth of the order book. Thus, we avoid the issues of selecting the risk-aversion parameters of the agents. We formulate the stochastic control problem of the exchange and prove a verification theorem characterizing the optimal incentives function. The main result of this article is that, given the form of the control problem of the exchange, the Backward Stochastic Partial Differential Equation (BSPDE) characterizing the problem boils down in fact to a linear parabolic Partial Differential Equation (PDE). This simplifies drastically the control problem, as it admits a classical Feynman-Kac representation. Moreover, with the choice of suitable intensity and penalty functions, the asymptotic solution of the control problem is obtained in closed-form by solving a simple second order linear ordinary differential equation (ODE) with constant coefficients. It is of particular interest for the practical design of incentives as they are now analytically obtained as functions of the distance from the mid-price and the market parameters. In particular, the flexibility of the model allows the exchange to modify specific limits of the order book and to design various forms of incentives (increase overall liquidity, more emphasis on first and second limit etc.). We illustrate numerically the form of the incentives produced by the model and its influence on the shape of the limit order book (i.e. on market liquidity) by simulating the associated SPDE controlled by the optimal incentives. We also study numerically the dependence of our model to the different market parameters.\\

The article is designed as follows. In Section \ref{sec_model}, we describe the order book dynamics and the optimization problem of the exchange. In particular, we propose a slightly modified version of the limit order book model introduced in \cite{cont2021stochastic}, to take into account the impact of incentives. In Section \ref{sec_mainresult}, we provide a BSPDE representation of the optimal control problem, and show that it admits a unique solution. We also prove a verification theorem characterizing the optimal control of the exchange. Section \ref{sec_longtermbehavior} is devoted to the asymptotic analysis of the control problem and the associated closed-form solutions, while Section \ref{sec_numerical} outlines the main numerical results of the incentives model. The proofs are relegated to the Appendix. 

\section{The model}\label{sec_model}

Throughout the paper, we consider a filtered probability space $\left( \Omega, \mathcal F, \mathbb P; \mathbb F = (\mathcal F_t)_{t \ge 0} \right)$ satisfying the usual conditions. We assume this probability space to be large enough to support all the processes we introduce.

\subsection{A Cont-M\"uller \cite{cont2021stochastic} limit order book}

Let us introduce a market consisting of a single asset, whose price process is denoted by $(S_t)_{t \in [0,T]}$ over a time horizon $[0,T]$, with $T>0.$ We consider a model for the limit order book (LOB) based on the work of Cont and M\"uller \cite{cont2021stochastic}. In this paper, the authors introduce a LOB modelled by a density $u:(t,x) \in [0,T] \times [-L,L] \mapsto u_t(x) $ representing the volume of orders available at time $t$ at price $S_t + x$ per unit price, at the bid for $x<0$ and at the ask for $x>0$. For mathematical convenience, we limit the range of the argument $x$ to a bounded interval $[-L,L]$, with $L>0$, meaning that one cannot propose a bid (resp. ask) price lower (resp. higher) than $S_t-L$ (resp. $S_t+L$).\footnote{As $L$ can be taken arbitrarily large, this assumption is not restrictive for practical applications.} The function $u$ satisfies the following stochastic partial differential equation (SPDE):\footnote{$u_t$ is the density of the order book at time $t$, i.e. the volume available at the first limit at the ask (for instance) at time $t$ is given by $\int_0^\delta u_t (x) dx$, where $\delta$ is the tick size of the asset. See \cite{cont2021stochastic} for more details.}

\begin{align}\label{SPDE0}
    du_t(x) = \left( \vphantom{\beta_b \nabla u_t(x) + \alpha_b u_t(x) -f^b(x)} \eta_a \Delta u_t(x) + \beta_a \nabla u_t(x) + \alpha_a u_t(x) +f^a(x) \right)dt + \sigma_a u_t(x) dW^a_t \quad \forall x \in (0,L),\\
    du_t(x) = \left( \eta_b \Delta u_t(x) - \beta_b \nabla u_t(x) + \alpha_b u_t(x) -f^b(x) \right)dt + \sigma_b u_t(x) dW^b_t \quad \forall x \in (-L,0),\nonumber
\end{align}
with boundary conditions:
\begin{align}\label{bound0}
    u_t(0^+) = u_t(0^-) = u_t(-L) = u_t(L) = 0 \quad \forall t \in [0,T],
\end{align}
with $u_0:[-L,L] \rightarrow u_0(x)$ fixed, and by convention $u_t(x) \le 0$ for $x<0$ and $u_t(x)\ge0$ for $x>0$, where, as mentioned above, $x$ is the (signed) distance from the mid-price. The process $(W^a_t, W^b_t)_{t \in [0,T]}$ is a 2-dimensional Brownian motion (possibly correlated), and the parameters are taken as follows : $\eta_a, \eta_b, \sigma_a, \sigma_b >0$, $\beta_a, \beta_b \ge 0$ and $\alpha_a, \alpha_b \le 0$. $f^a:(0,L) \rightarrow \mathbb [0,+\infty)$ and $f^b:(-L,0) \rightarrow [0,+\infty)$ are two functions. The interpretation is given below:
\begin{itemize}
    \item $f^b$ (resp. $f^a$) represents the rate of buy (resp. sell) order submissions at a distance $x$ from the mid-price;
    \item the term $\alpha_b u_t(x)$ (resp. $\alpha_b u_t(x)$) represents the proportional cancellation of limit buy (resp. sell) orders at a distance $x$ from the mid-price;
    \item the convection term $- \beta_b \nabla u_t(x)$ (resp. $\beta_a \nabla u_t(x)$) represents the replacement of buy (resp. sell) orders by orders closer to the mid-price;
    \item the diffusion term $\eta_b \Delta u_t(x)$ (resp. $\eta_a \Delta u_t(x)$) represents the cancellation and symmetric replacement of orders at a distance $x$ from the mid-price;
    \item the multiplicative noise term $\sigma_b u_t(x) dW^b_t$ (resp. $\sigma_a u_t(x) dW^a_t$) accounts for the high-frequency submissions and cancellations associated with HFT orders.
\end{itemize}

We will now considered a slightly modified version of the above model, to take into account the use of incentives from the exchange to optimize the shape of the LOB. More precisely, we assume that $u$ is now solution to the following SPDE:

\begin{align}\label{SPDE1}
    du_t(x) = \left( \eta_a \Delta u_t(x) + \beta_a \nabla u_t(x) + \alpha_a u_t(x) +f^a(x, Z^a_t(x)) \right)dt + \sigma_a u_t(x) dW^a_t \ \forall x \in (0,L),\\
    du_t(x) = \left( \eta_b \Delta u_t(x) - \beta_b \nabla u_t(x) + \alpha_b u_t(x) -f^b(x, Z^b_t(x)) \right)dt + \sigma_b u_t(x) dW^b_t \ \forall x \in (-L,0),\nonumber
\end{align}
where $f^a, f^b: (0,L) \times [0,+\infty) \rightarrow [0,+\infty)$ verify:
\begin{itemize}
    \item $f^a$ and $f^b$ are $C^1$;
    \item $f^a(x,.)$ and $f^b(x,.)$ are increasing for all $x \in (0,L)$;
    \item $f^a(x,.)$ and $f^b(x,.)$ are concave for all $x \in (0,L)$.
\end{itemize}

\begin{rem}
Although $f^a$ and $f^b$ are typically decreasing in $x$, our model does not require it. 
\end{rem}

In Equation \eqref{SPDE1}, the processes $Z^a_t(x)$ et $Z^b_t(x)$ are incentives that the exchange commits to give to the trader placing a unit order at a distance $x$ from the mid-price at time $t$, if the order is executed (hence the incentive is not received by the trader at time $t$ but at the time of the trade). They belong respectively to the sets of admissible controls $\mathcal A^a$ and $\mathcal A^b$ defined by
\begin{align*}
\mathcal A^j\! :=\! \bigg\{ Z:\Omega \times [0,T] \times (0,L) \rightarrow [0,+\infty) \Big| Z &\text{ is } \mathcal P \otimes \mathcal B \big( (0,L) \big)\text{-measurable,}\\
&\text{ and } \mathbb E\! \left[\!\int_0^T \!\!\! \int_0^L f^j \left(x, Z_t(x) \right)^2 dx dt \right]\! <\!+\infty \! \bigg\}    
\end{align*}

for $j \in \{a,b\},$ where $\mathcal P$ denotes the $\sigma$-algebra of $\mathbb F$-predictable subsets of $\Omega \times [0,T]$, and $\mathcal B \big( (0,L) \big)$ denotes the Borelian sets of $(0,L)$.

\subsection{The optimization problem}

The exchange has two symmetrical problems, on the bid side and ask side of the LOB, respectively. We only consider here the problem on the ask side (the problem on the bid side is treated similarly). The exchange aims at maximizing
\begin{align*}
    J(Z^a) := \mathbb E \left[ \int_0^T \left( \int_0^L \left(u^{Z^a}_t(x) - g(x,Z^a_t(x)) \right) dx \right)dt \right],
\end{align*}
over the set of admissible controls $\mathcal A^a$, where $g:(0,L) \times [0,+\infty) \rightarrow [0,+\infty)$ is a penalty function which verifies the following assumptions:
\begin{itemize}
    \item $g(x,0) = 0$ for all $x \in (0,L)$;
    \item $g$ is nondecreasing in both its variables;
    \item $g$ is $C^1$;
    \item $g(x,.)$ is convex for all $x \in (0,L)$.
\end{itemize}

Hence, the exchange wants to find $\hat Z^a \in \mathcal A^a$ such that
\begin{align}\label{askPB}
J(\hat Z^a) = \underset{Z^a \in \mathcal A^a}{\sup} J(Z^a). 
\end{align}

\begin{rem}
Ideally, the exchange would like the LOB to be as full as possible. However, providing incentives induces a cost for the exchange. The function $g$ is a penalty function that allows to take this cost into account. The first assumption $g(x,0) = 0$ is perfectly natural,  as providing no incentive induces no cost. Similarly, $g$ should be nondecreasing in its second variable, as a higher incentive induces a higher cost. We also assume that $g$ is nondecreasing in its first variable, to model the fact that the exchange is mainly interested in having a lot of liquidity at the first limits of the LOB, and is therefore ready to offer a higher incentive at those limits -- basically, we penalize ‘‘less’’ the cost induced by incentives at the first limits than the cost induced by incentives at limits further away from the mid-price.
\end{rem}

\section{Main results}\label{sec_mainresult}

Let us define the Hamiltonian $H:(0,L) \times \mathbb R \times [0,+\infty) \times \mathbb R \times \mathbb R \rightarrow \mathbb R$ by
$$H(x,u,z,p,q) = u - g(x,z) + (\alpha_a u + f^a(x,z) )p + \sigma_a u q.$$
From the properties of $g$ and $f^a$, we get the following lemma:
\begin{lemme}\label{Hconcave}
For all $(x,p,q) \in [0,L] \times [0,+\infty) \times \mathbb R$, the function
$$(u,z) \in [0,+\infty) \times [0,+\infty) \mapsto H(x,u,z,p,q)$$
is concave.
\end{lemme}
We introduce the differential operator $\mathcal L$ given by
$$\mathcal L u_t(x) = \eta_a \Delta u_t(x) + \beta_a \nabla u_t(x).$$
Let us then consider the following backward stochastic partial differential equation (BSPDE) $p(t,x) \in \mathbb R$ and $q(t,x) \in \mathbb R$:
\begin{align}\label{BSPDE0}
    dp_t(x) = - \left\{ \mathcal L^* p_t(x) + \partial_u H\left(x,u_t(x), Z^a_t(x), p_t(x), q_t(x) \right) \right\} dt + q_t(x) dW^a_t, \quad \forall (t,x) \in (0,T) \times (0,L),
\end{align}
where $\mathcal L^*$ is the adjoint operator of $\mathcal L$, with boundary conditions:
\begin{align}\label{boundBSPDE}
    p_T(x) = 0, \quad \forall x \in (0,L),\\
    p_t(x) = 0, \quad \forall (t,x) \in (0,T) \times \{0,L\}.\nonumber
\end{align}
Note that we can write the BSPDE \eqref{BSPDE0} as
\begin{align}\label{BSPDE1}
    dp_t(x) = - \left\{ \eta_a \Delta p_t(x) - \beta_a \nabla p_t(x) + \left(1 + \alpha_a p_t(x)+ \sigma_a q_t(x) \right) \right\} dt + q_t(x) dW^a_t, \quad \forall (t,x) \in (0,T) \times (0,L).
\end{align}
The following result is proved in Appendix \ref{Proof1}:
\begin{thm}\label{maximum}
Let $\hat Z ^a \in \mathcal A^a$, and denote by $\hat u$ the associated solution to SPDE \eqref{SPDE1}. Assume there exists a solution $(\hat p, \hat q)$ to the BSPDE  \eqref{BSPDE1} with boundary conditions \eqref{boundBSPDE}. Assume also that:
\begin{enumerate}[label=(\roman*)]
    \item For all $t \in [0,T)$ and $x \in (0,L)$, 
    $$H(x,\hat u_t(x), \hat Z^a_t(x), \hat p_t(x), \hat q_t(x)) = \underset{z \in [0,+\infty)}{\sup}\ H(x,\hat u_t(x), z, \hat p_t(x), \hat q_t(x));$$
    \item For all $t \in [0,T]$ and $x \in [0,L]$,
    $$\hat p_t(x) \ge 0;$$
    \item For all $Z^a \in \mathcal A^a$,
    $$\mathbb E \left[ \int_0^L \int_0^T \left(u^{Z^a}_t(x) - \hat u_t(x) \right)^2 \hat q_t(x)^2 dt dx \right]<+\infty,$$
    where $u^{Z^a}$ is the solution to \eqref{SPDE1} controlled by $Z^a$;
    \item For all $Z^a \in \mathcal A^a$,
    $$\mathbb E \left[ \int_0^L \int_0^T \left( u^{Z^a}_t(x)\right)^2 \hat p_t(x)^2 dt dx \right]<+\infty;$$
\end{enumerate}
then $\hat Z^a$ is an optimal control for the control problem \eqref{askPB}.
\end{thm}

As the boundary conditions \eqref{boundBSPDE} are deterministic, we look for a solution to the BSPDE \eqref{BSPDE1} such that $q(t,x) = 0$ for all $(t,x) \in [0,T] \times [0,L]$. In that case, $p$ solves the following linear parabolic PDE:\footnote{From now on, we write $p(t,x)$ instead of $p_t(x)$ to highlight the fact that $p$ is deterministic.}
\begin{align}\label{PDE0}
    \partial_t p(t,x) + \eta_a \Delta p(t,x) - \beta_a \nabla p(t,x)  + \alpha_a p(t,x) + 1 = 0,
\end{align}
with boundary conditions:
\begin{align}\label{boundPDE0}
    p(T,x) = 0 \ \forall x \in (0,L)\\
    p(t,x) = 0 \ \forall (t,x) \in (0,T) \times \{0,L\}. \nonumber 
\end{align}
This is a classical boundary value problem, and the solution admits a Feynman-Kac representation of the following form:
\begin{align} \label{FK0}
p(t,x) = \mathbb E \left[ \int_t^{T\wedge \gamma^{t,x}} e^{\alpha_a (s-t)}ds\right],    
\end{align}
where we consider a new Brownian motion $\widehat W$ and the following SDE:
\begin{align}\label{SDEX}
dX_s = -\beta_a ds + \sqrt{2 \eta_a}d\widehat W_s,  
\end{align}
and $\gamma^{t,x}$ denotes the stopping time
$$\gamma^{t,x} := \inf \left\{ s\ge t \left| X^{t,x}_s \not\in (0,L) \right.\right\},$$
where $(X^{t,x}_s)_{s\ge 0}$ is the solution to SDE \eqref{SDEX} starting at time $t$ and at point $x$.\\

As $H$ is concave in $z$, the optimal control is then given by the solution to:
\begin{align*}
    & \partial_z H\left(x,u_t(x), z, p(t,x), q(t,x) \right)= 0\\
     \iff & - \partial_z g(x,z) + p(t,x) \partial_z f^a(x,z)=0,
\end{align*}
if such a solution exists.

\begin{rem}
If the above equation has a solution for all $t$ and $x$, it is clear that the corresponding optimal control satisfies condition $(i)$ in Theorem \ref{maximum}. Furthermore, notice from Equation \eqref{FK0} that $p(t,x) \ge 0$ for all $t$ and $x$, hence condition $(ii)$ is also satisfied. Note that we obtain here a solution such that $q(t,x) = 0$ for all $t$ and $x$, and therefore condition $(iii)$ is also naturally verified. Finally, as $p$ is clearly bounded from Equation \eqref{FK0}, Theorem 6.7 in \cite{da2014stochastic} guarantees that condition $(iv)$ is satisfied.
\end{rem}

\section{Long-term behaviour and closed-form formula}\label{sec_longtermbehavior}

\subsection{Asymptotics of $p$}

In order to study the asymptotic behaviour of the function $p$, we consider the stationary equation:
\begin{align}\label{PDELT}
    \eta_a \Delta \bar p(x) - \beta_a \nabla \bar p(x)  + \alpha_a \bar p(x) + 1 = 0,
\end{align}
with boundary conditions:
\begin{align}\label{boundPDELT}
    \bar p(x) = 0 \ \forall x \in (0,L).
\end{align}
The Feynman-Kac representation of the solution is given by:
\begin{align} \label{FKLT}
\bar p(x) = \mathbb E \left[ \int_0^{\gamma^{0,x}} e^{\alpha_a s}ds\right].  
\end{align}

The following proposition is straightforward by dominated convergence:

\begin{prop}\label{asymP}
For all $x \in [0,L]$, we have 
$$p(t,x) \underset{T\rightarrow + \infty}{\longrightarrow} \bar p(x).$$
\end{prop}

Moreover, notice that Equation \eqref{PDELT} is a simple second order linear ODE, and together with conditions \eqref{boundPDELT}, its solution is given explicitly by:
\begin{align}\label{LTP}
    \bar p(x) = -\frac 1{\alpha_a} + \mu_+ e^{\nu_+ x} + \mu_- e^{\nu_- x}
\end{align}
for all $x \in [0,L]$, with 
$$\nu_+ = \frac{\beta_a + \sqrt{\beta_a^2 - 4\eta_a \alpha_a}}{2\eta_a} \qquad \text{and} \qquad \nu_- = \frac{\beta_a - \sqrt{\beta_a^2 - 4\eta_a \alpha_a}}{2\eta_a},$$
and
$$\mu_+ = \frac 1{\alpha_a} - \mu_- \qquad \text{and} \qquad \mu_- = \frac{1}{\alpha_a} \frac{e^{\nu_+L} - 1}{e^{\nu_+ L} - e^{\nu_- L}}.$$

\subsection{An example with closed-form solution}

We propose a penalty function $g$ and an intensity function $f^a$ with the following forms:
$$g(x,z) = \bar A ze^{\bar a x},$$
and
$$f^a (x,z) =  \lambda z^r e^{-\kappa x} + \lambda_0(x)$$
for all $(x,z) \in  (0,L) \times [0,+\infty)$, where $\bar A, \bar a, \lambda, \kappa >0$, $r \in (0,1)$, and $\lambda_0$ is a non-negative function representing the rate of sell order submissions at a distance $x$ from the mid-price in the absence of incentives.\\

\begin{rem}
These two functions deserve a few comments:
\begin{itemize}
    \item For a fixed $x \in (0,L)$, the intensity function $f^a(x,.)$ has both a constant part, $\lambda_0(x)$, corresponding to the intensity of orders arrival in the absence of incentives (when $\left(Z^a_t(x)\right)_{t\in [0,T]}$ is identically equal to $0$), and a part that depends on the incentive: $\lambda z^r e^{-\kappa x}$.This second part is exponentially decreasing in $x$ to model the fact that the incentive has less and less impact as we move away from the mid-price: indeed, the incentive is only perceived by the liquidity provider if her order is executed, which is more likely when the order is close to the mid-price. It is concave in $z$ to model the decrease of the marginal liquidity gain derived from each additional unit of incentive.
    \item It seems natural to choose a penalty function $g$ that is linear in $z$, as $g$ should intuitively be proportional to the cost incurred for the platform. It is increasing in $x$ to model the fact that the exchange is more inclined to pay to increase liquidity at the first limits than far away from the mid-price. Here, we take the example of an exponential increase, but one could reasonably choose to work with any function $g(x,z) = z h(x)$ for any $C^1$ increasing function $h$ (the choice of $h$ depends on how much the platform wishes to concentrate liquidity on the first limits).
\end{itemize}
\end{rem}

\begin{prop}
In this framework, the optimal control $\hat Z^a \in \mathcal A^a$ is given by
$$\hat Z^a_t(x) = \left( \frac{p(t,x) \lambda r}{\bar A e^{(\bar a + \kappa )x}} \right)^{\frac 1{1-r}}.$$
When $T \rightarrow +\infty$, this control reaches a stationary state given in closed-form by
\begin{align}\label{closed}
   \hat Z^a(x) = \left( \frac{\bar p(x) \lambda r}{\bar A e^{(\bar a + \kappa )x}} \right)^{\frac 1{1-r}}.
\end{align}
\end{prop}

\section{Numerical simulations}\label{sec_numerical}
\subsection{Model parameters}

In this section, we apply our model to the case of a (fictitious) symmetric limit order book with a tick size given by $\delta = 0.01\ \$$, with the following parameters:
\begin{itemize}
    \item Diffusion parameters: $\eta_b = \eta_a = 10^{-3}\ \$^2 \cdot \min^{-1}$.
    \item Convection parameters: $\beta_b = \beta_a = 2 \cdot 10^{-2}\ \$ \cdot \min^{-1}$.
    \item Cancellation parameters: $\alpha_b = \alpha_a = -0.2 \min^{-1}$.
    \item Noise parameters: $\sigma_b = \sigma_a = 0.3 \min^{-\frac 12}$.
    \item Correlation between $W^a$ and $W^b$: $\rho = -0.05$.
    \item Rate of submissions of the form
    $$f^a(x,z) = f^b(-x,z) = \lambda \sqrt{z} e^{-\kappa x } + \lambda_0 e^{-\kappa_0 x}, $$
    with:
    \begin{itemize}
        \item[] $\lambda_0 = 50000\ \$^{-1} \cdot \min^{-1},$
        \item[] $\kappa_0 = 50\ \$^{-1},$
        \item[] $\lambda = 630000\ \$^{-\frac 32} \cdot \min^{-1},$
        \item[] $\kappa = 100\ \$^{-1}.$
    \end{itemize}
    This means that, in the absence of incentives, unit limit orders arrive at a distance of $0.01\ \$$ from the mid-price at a constant rate of $\delta \lambda_0 e^{-0.01 \kappa_0 } \simeq 303 \min^{-1}, $ and at a distance of $0.05\ \$$ from the mid-price at a constant rate of $\delta \lambda_0 e^{-0.05 \kappa_0 } = 41 \min^{-1}.$ An incentive of $z=0.01\ \$$ per unit order increases the rate at a distance of $0.01\ \$$ from the mid-price by $\delta \lambda \sqrt{z} e^{-0.01 \kappa } \simeq 232 \min^{-1}, $ and at a distance of $0.05\ \$$ from the mid-price by $\delta \lambda \sqrt{z} e^{-0.05 \kappa } = 4 \min^{-1}.$
\end{itemize}
Regarding the objective function we consider the following:
\begin{itemize}
    \item Boundary on the LOB: $L=0.11 \ \$$ (this corresponds to $10$ limits on each side).
    \item Time horizon: $T=30 \min$. This horizon ensures convergence toward the stationary incentives at time $t=0$.
    \item Penalty function of the form $g(x,z) = \bar A ze^{\bar a x},$ with:
        \begin{itemize}
        \item[] $\bar A = 4200\ \$^{-2},$
        \item[] $\bar a = 50\ \$^{-1}.$
    \end{itemize}
\end{itemize}

\begin{rem}
The parameters $\eta_b, \eta_a, \beta_b, \beta_a, \alpha_b, \alpha_a, \sigma_b, \sigma_a, \lambda_0,$ and $\kappa_0$ can easily be calibrated on historical data. The parameters $\lambda$ and $\kappa$ can be calibrated if the exchange already has historical data with different levels of incentives; otherwise, some tests and heuristics should be used. Finally, $\bar A$ and $\bar a$ should be calibrated depending on the desired shape for the limit order book and the amount the exchange is ready to spend; simulations can be carried out in order to find the appropriate values for these parameters. 
\end{rem}

\subsection{Optimal incentives and shape of the order book}

\begin{figure}[!h]\centering
\includegraphics[width=0.87\textwidth]{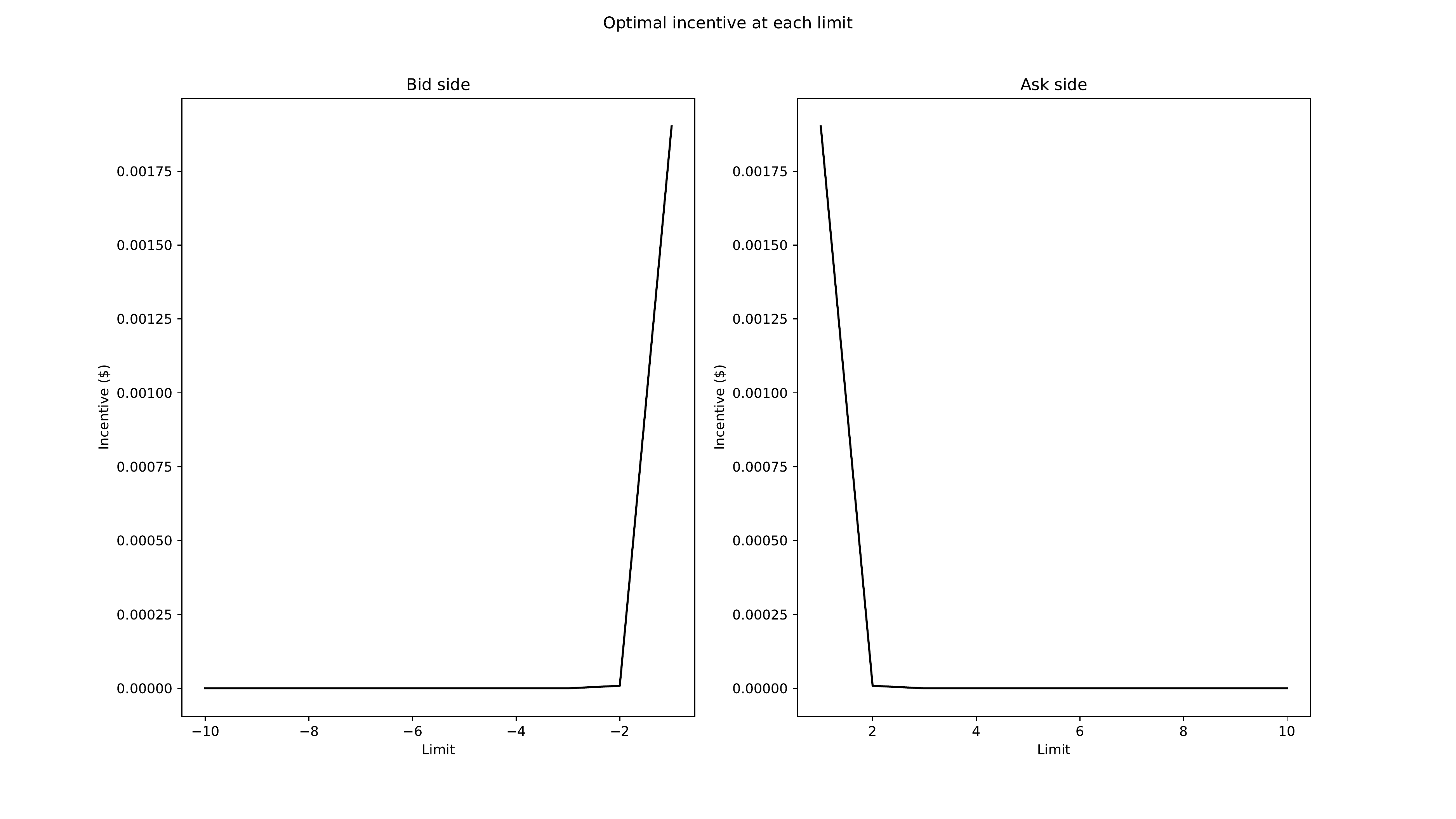}\\
\caption{Optimal incentive for a unit order at each limit of the limit order book.}\label{opt_0}
\end{figure}

We compute the asymptotic optimal incentives using Equation \eqref{closed}, and plot them in Figure \ref{opt_0}. As we observe, the optimal incentives decrease at a very fast rate as we get away from the mid-price. For more clarity, we report these incentives in Table \ref{table_0} (the incentives are symmetric at the bid and at the ask). This fast decrease is expected, as the exchange wants to concentrate liquidity at the first limits of the order book.\\

\begin{figure}[!h]\centering
\includegraphics[width=0.88\textwidth]{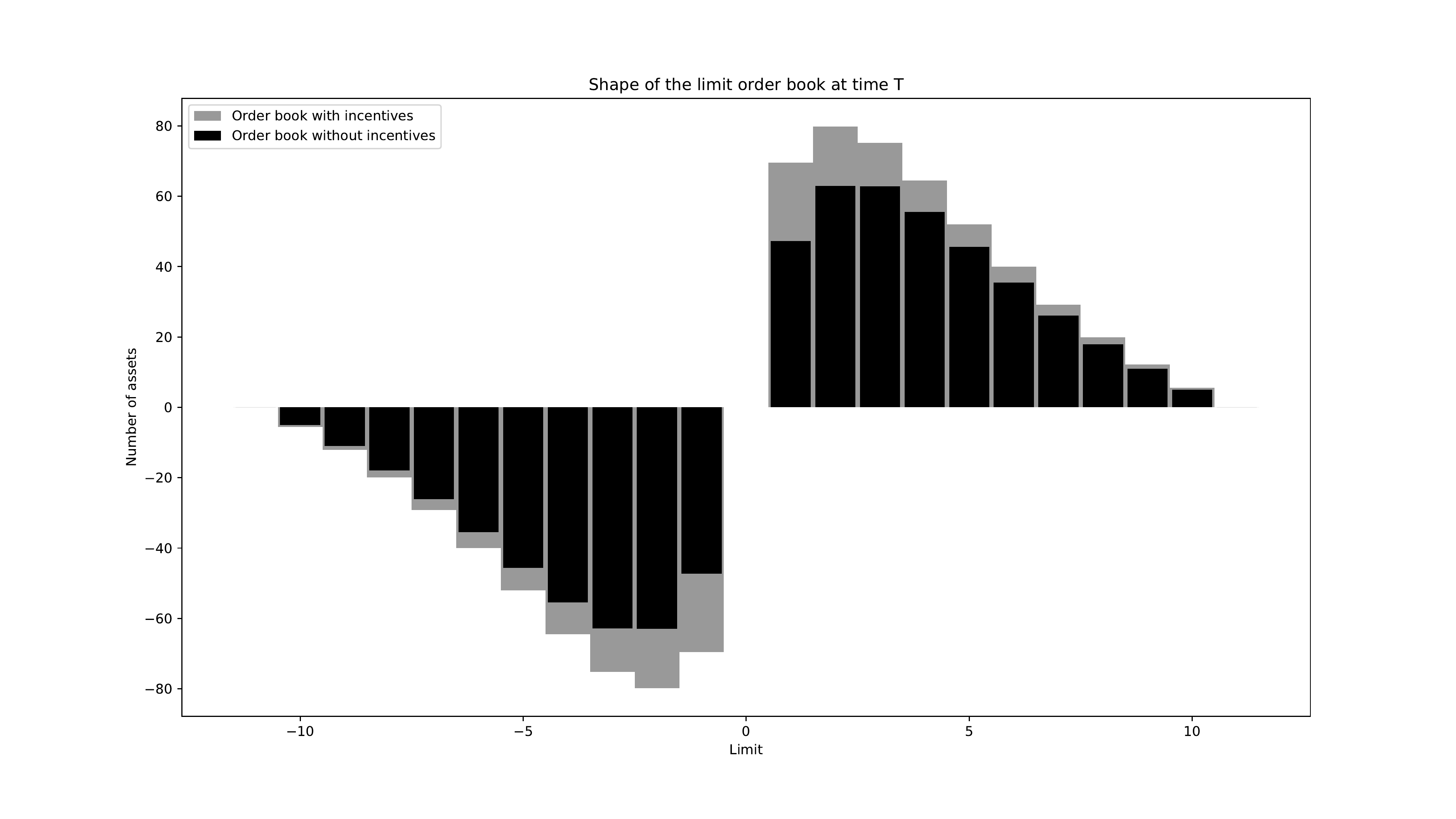}\\
\caption{Shape of the limit order book at time $T$.}\label{shape_0}
\end{figure}

\begin{table}[!h]
\begin{center}
\begin{tabular}{c  c} 
 \hline 
Limit & Incentive (\$) \\ [0.5ex] 
 \hline
 $1st$ & $1.90 \cdot 10^{-3}$ \\ [0.5ex] 
 $2nd$ &  $8.56 \cdot 10^{-6}$ \\  [0.5ex] 
 $3rd$ &  $2.11 \cdot 10^{-8}$ \\ [0.5ex] 
 $4th$ &  $4.01 \cdot 10^{-11}$ \\ [0.5ex] 
 $5th$ &  $6.41 \cdot 10^{-14}$ \\ [0.5ex] 
 $6th$ &  $8.96 \cdot 10^{-17}$ \\ [0.5ex] 
 $7th$ &  $1.09 \cdot 10^{-19}$ \\ [0.5ex] 
 $8th$ &  $1.12 \cdot 10^{-22}$ \\ [0.5ex] 
 $9th$ &  $8.86 \cdot 10^{-26}$ \\ [0.5ex] 
 $10th$ &  $3.84 \cdot 10^{-29}$ \\ [0.5ex] 
 \hline 
\end{tabular}
\end{center}
\caption {Optimal incentive for a unit order at each limit of the limit order book.}
\label{table_0}
\end{table}
\newpage
We then simulate the limit order book given by equation \eqref{SPDE1} with and without incentives, to obtain its average shape at time $T$ and compare the two cases. The shape is plotted in Figure \ref{shape_0}. We see that the use of incentives allows to increase the liquidity level of the order book, particularly at the first limits. Note that, in the absence of incentives, the first limit is only the $5$-th more filled of the order book, whereas with incentives it becomes the $3$-rd more filled. 

\subsection{Impact of the different parameters}

We now play with the parameters of Equation \eqref{SPDE1} to study their impact on the optimal incentives and on the limit order book.

\subsubsection{Diffusion parameter}

\begin{figure}[!h]\centering
\includegraphics[width=\textwidth]{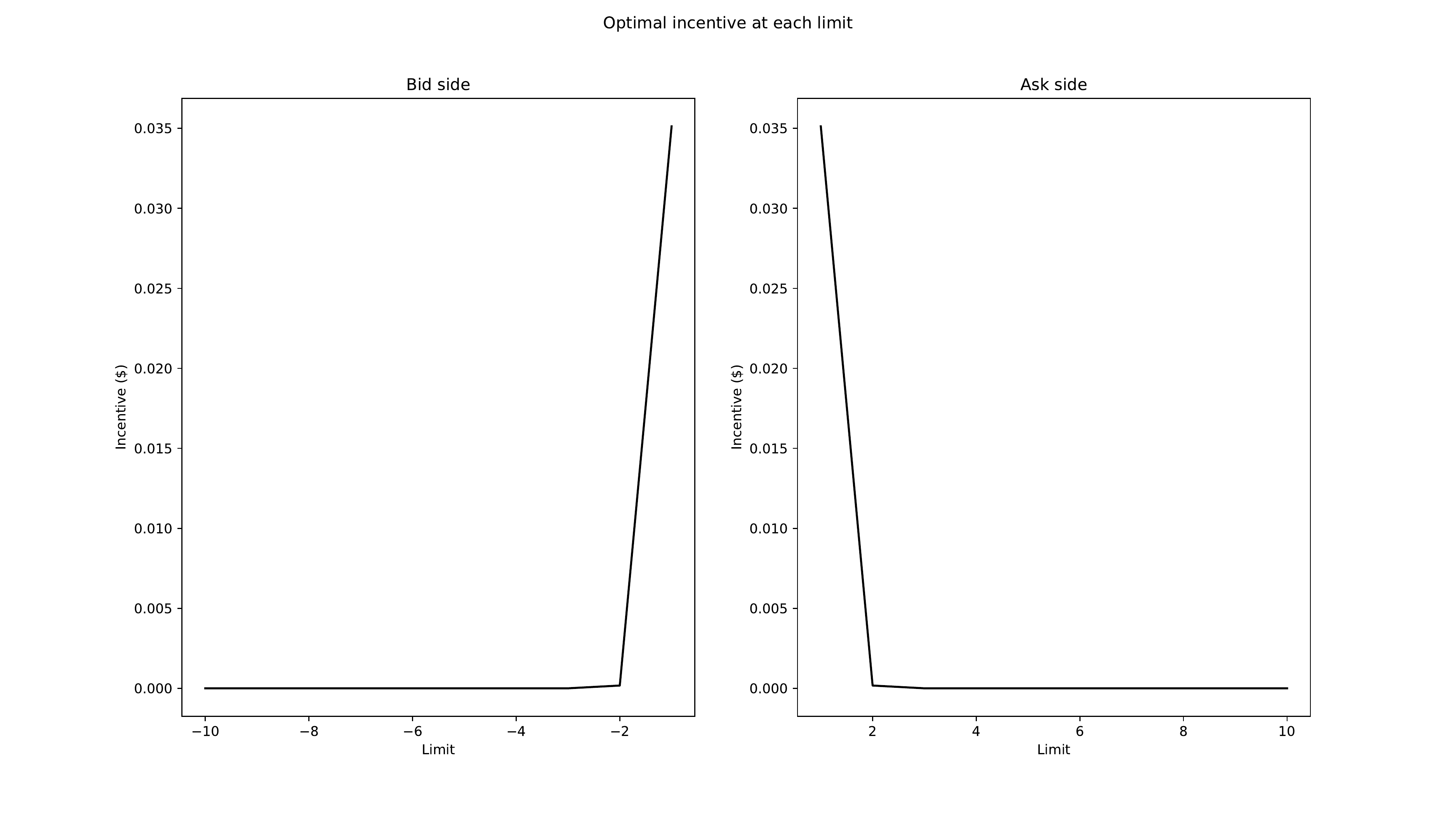}\\
\caption{Optimal incentive for a unit order at each limit of the limit order book.}\label{opt_1}
\end{figure}

Let us start with the diffusion parameter. We now divide it by 2, and set $\eta_b = \eta_a = 5 \cdot 10^{-4}\ \$^2 \cdot \min^{-1}$. The optimal incentives are now plotted in Figure \ref{opt_1} and reported in Table \ref{table_1}.\\

\begin{figure}[!h]\centering
\includegraphics[width=\textwidth]{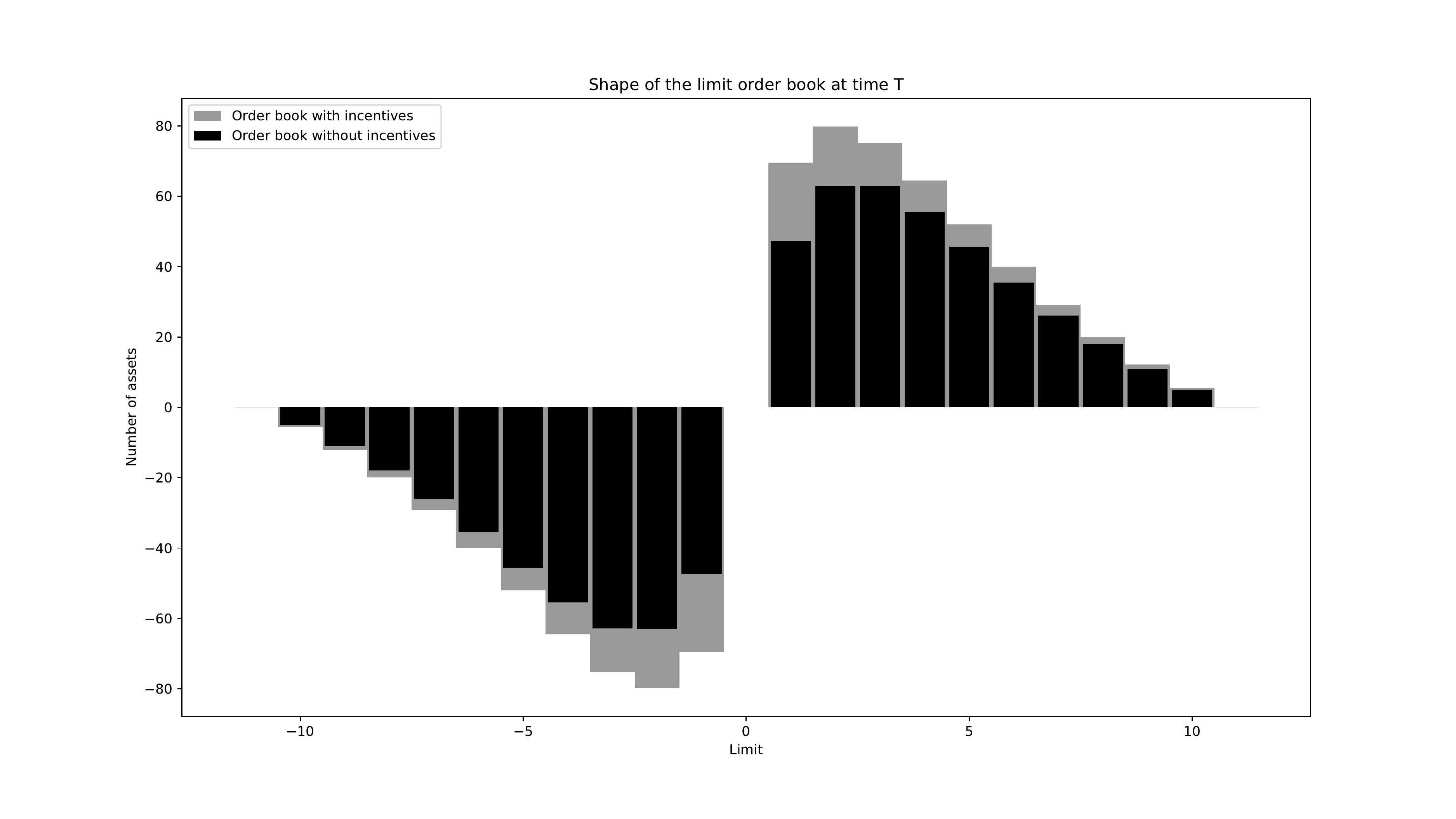}\\
\caption{Shape of the limit order book at time $T$.}\label{shape_1}
\end{figure}

We observe that a decrease in the diffusion parameter leads to an increase in the incentives. Indeed, reducing $\eta_b$ and $\eta_a$ improves the impact of the incentives, as the new orders posted are less likely to be replaced by orders further away from the mid-price. This is confirmed by the shape of the order book in Figure \ref{shape_1}, where we see that the overall level of the order book is higher than in Figure \ref{shape_0}, particularly at the first limits. 

\newpage
\begin{table}[!h]
\begin{center}
\begin{tabular}{c  c} 
 \hline 
Limit & Incentive (\$) \\ [0.5ex] 
 \hline
 $1st$ & $3.51 \cdot 10^{-2}$ \\ [0.5ex] 
 $2nd$ &  $1.73 \cdot 10^{-4}$ \\  [0.5ex] 
 $3rd$ &  $4.77 \cdot 10^{-7}$ \\ [0.5ex] 
 $4th$ &  $1.03 \cdot 10^{-9}$ \\ [0.5ex] 
 $5th$ &  $1.90 \cdot 10^{-12}$ \\ [0.5ex] 
 $6th$ &  $3.15 \cdot 10^{-15}$ \\ [0.5ex] 
 $7th$ &  $4.68 \cdot 10^{-18}$ \\ [0.5ex] 
 $8th$ &  $6.02 \cdot 10^{-21}$ \\ [0.5ex] 
 $9th$ &  $6.14 \cdot 10^{-24}$ \\ [0.5ex] 
 $10th$ &  $3.55\cdot 10^{-27}$ \\ [0.5ex] 
 \hline 
\end{tabular}
\end{center}
\caption {Optimal incentive for a unit order at each limit of the limit order book.}
\label{table_1}
\end{table}

\newpage

\subsubsection{Convection rate}

\begin{figure}[!h]\centering
\includegraphics[width=0.92\textwidth]{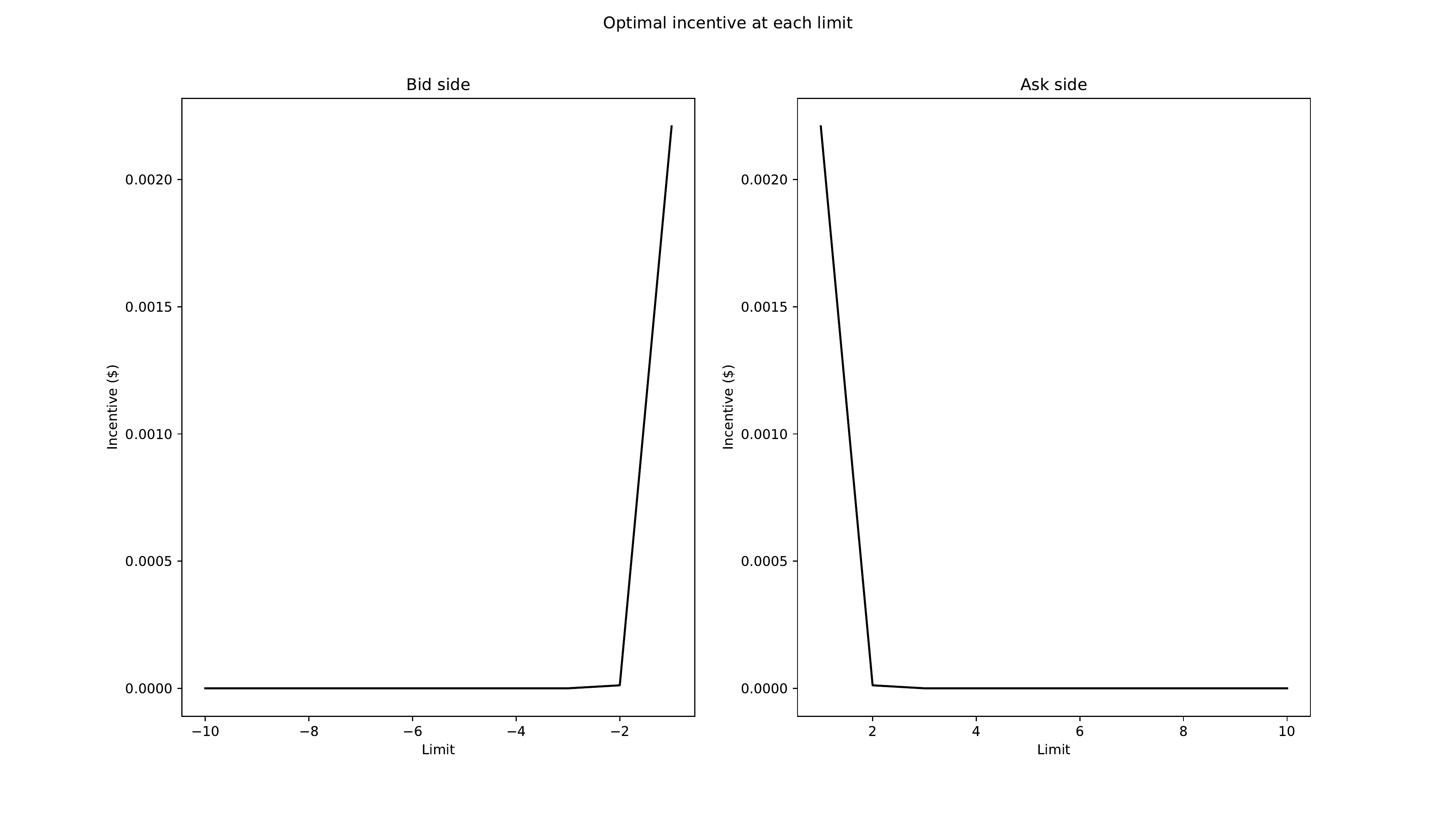}\\
\caption{Optimal incentive for a unit order at each limit of the limit order book.}\label{opt_2}
\end{figure}

We now turn to the convection rate and multiply it by 5, setting $\beta_b = \beta_a = 0.1\ \$ \cdot \min^{-1}$. The optimal incentives are now plotted in Figure \ref{opt_2} and reported in Table \ref{table_2}.\\

\begin{table}[!h]
\begin{center}
\begin{tabular}{c  c} 
 \hline 
Limit & Incentive (\$) \\ [0.5ex] 
 \hline
 $1st$ & $2.21 \cdot 10^{-3}$ \\ [0.5ex] 
 $2nd$ &  $1.18 \cdot 10^{-5}$ \\  [0.5ex] 
 $3rd$ &  $3.53 \cdot 10^{-8}$ \\ [0.5ex] 
 $4th$ &  $8.37 \cdot 10^{-11}$ \\ [0.5ex] 
 $5th$ &  $1.74 \cdot 10^{-13}$ \\ [0.5ex] 
 $6th$ &  $3.29 \cdot 10^{-16}$ \\ [0.5ex] 
 $7th$ &  $5.80 \cdot 10^{-19}$ \\ [0.5ex] 
 $8th$ &  $9.36 \cdot 10^{-22}$ \\ [0.5ex] 
 $9th$ &  $1.28 \cdot 10^{-24}$ \\ [0.5ex] 
 $10th$ &  $1.09 \cdot 10^{-27}$ \\ [0.5ex] 
 \hline 
\end{tabular}
\end{center}
\caption {Optimal incentive for a unit order at each limit of the limit order book.}
\label{table_2}
\end{table}

\begin{figure}[!h]\centering
\includegraphics[width=0.88\textwidth]{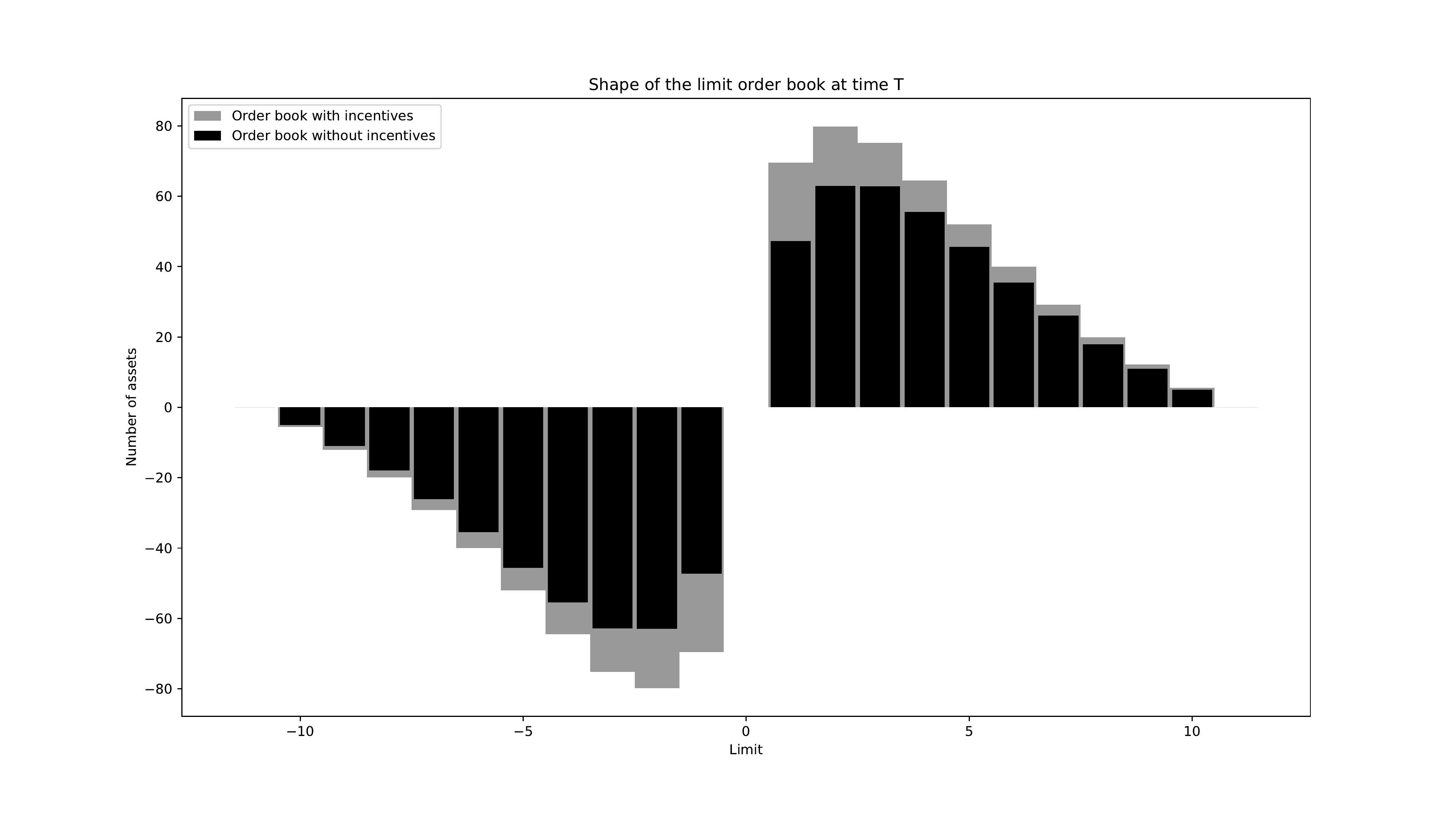}\\
\caption{Shape of the limit order book at time $T$.}\label{shape_2}
\end{figure}

We observe that an increase in the convection rate leads to an increase in the incentives. Indeed, increasing $\beta_b$ and $\beta_a$ means that orders will naturally move towards the first limit, hence increasing the incentives even at the last limits will have a good impact on liquidity. The overall level of the order book in Figure \ref{shape_2} is lower than the one in Figure \ref{shape_0}, because orders at the first limit are also more rapidly matched by a market order.

\newpage

\subsubsection{Cancellation rate}

\begin{figure}[!h]\centering
\includegraphics[width=0.88\textwidth]{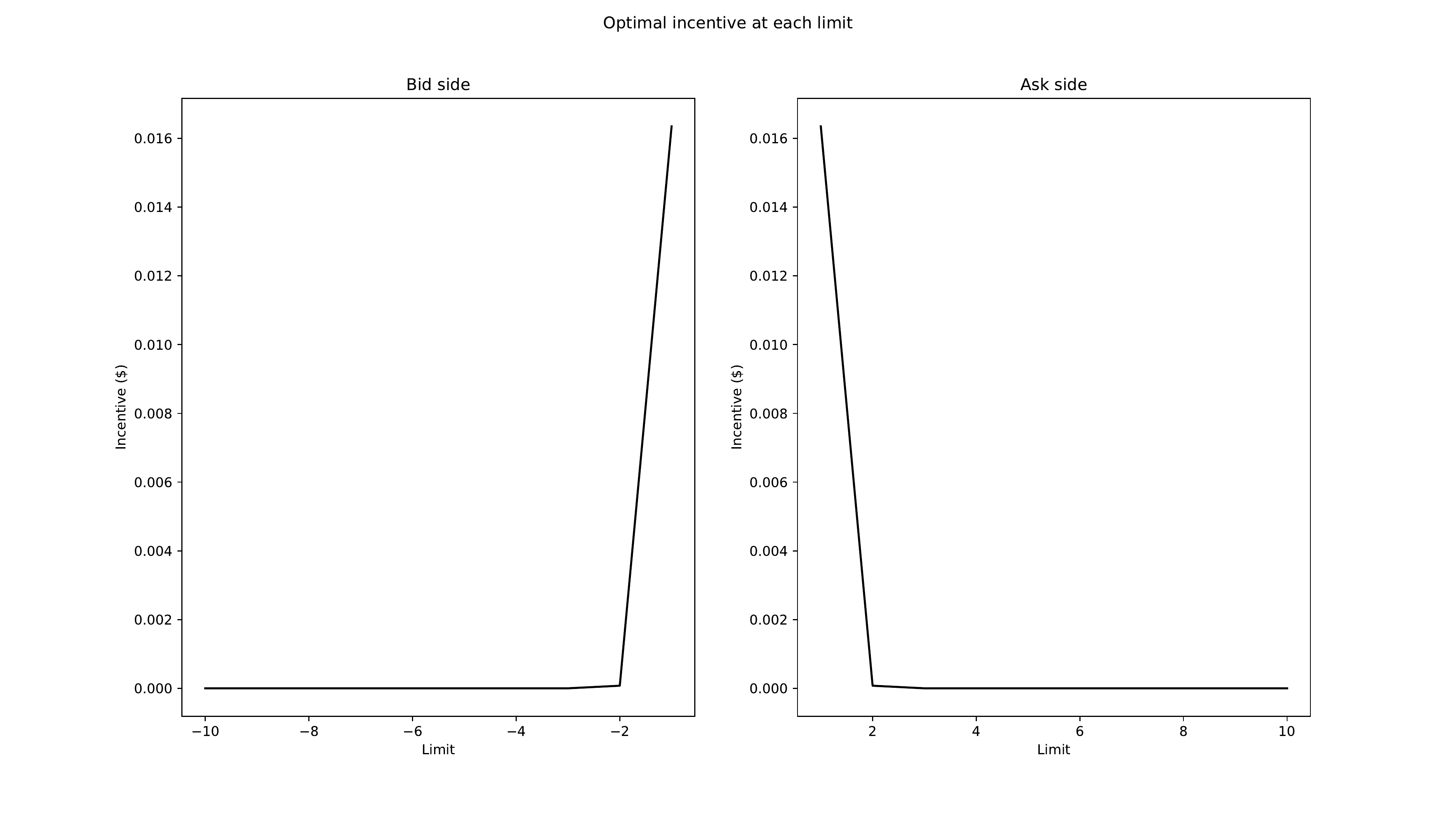}\\
\caption{Optimal incentive for a unit order at each limit of the limit order book.}\label{opt_3}
\end{figure}

Finally, we study the impact of the cancellation rate. We multiply it by 2, setting $\alpha_b = \alpha_a = -0.4 \min^{-1}$. The optimal incentives are now plotted in Figure \ref{opt_3} and reported in Table \ref{table_3}.\\

\begin{table}[!h]
\begin{center}
\begin{tabular}{c  c} 
 \hline 
Limit & Incentive (\$) \\ [0.5ex] 
 \hline
 $1st$ & $1.64 \cdot 10^{-2}$ \\ [0.5ex] 
 $2nd$ &  $7.57 \cdot 10^{-5}$ \\  [0.5ex] 
 $3rd$ &  $1.95 \cdot 10^{-7}$ \\ [0.5ex] 
 $4th$ &  $3.90 \cdot 10^{-10}$ \\ [0.5ex] 
 $5th$ &  $6.69 \cdot 10^{-13}$ \\ [0.5ex] 
 $6th$ &  $1.01 \cdot 10^{-15}$ \\ [0.5ex] 
 $7th$ &  $1.36 \cdot 10^{-18}$ \\ [0.5ex] 
 $8th$ &  $1.57 \cdot 10^{-21}$ \\ [0.5ex] 
 $9th$ &  $1.41 \cdot 10^{-24}$ \\ [0.5ex] 
 $10th$ &  $7.07 \cdot 10^{-28}$ \\ [0.5ex] 
 \hline 
\end{tabular}
\end{center}
\caption {Optimal incentive for a unit order at each limit of the limit order book.}
\label{table_3}
\end{table}

\begin{figure}[!h]\centering
\includegraphics[width=\textwidth]{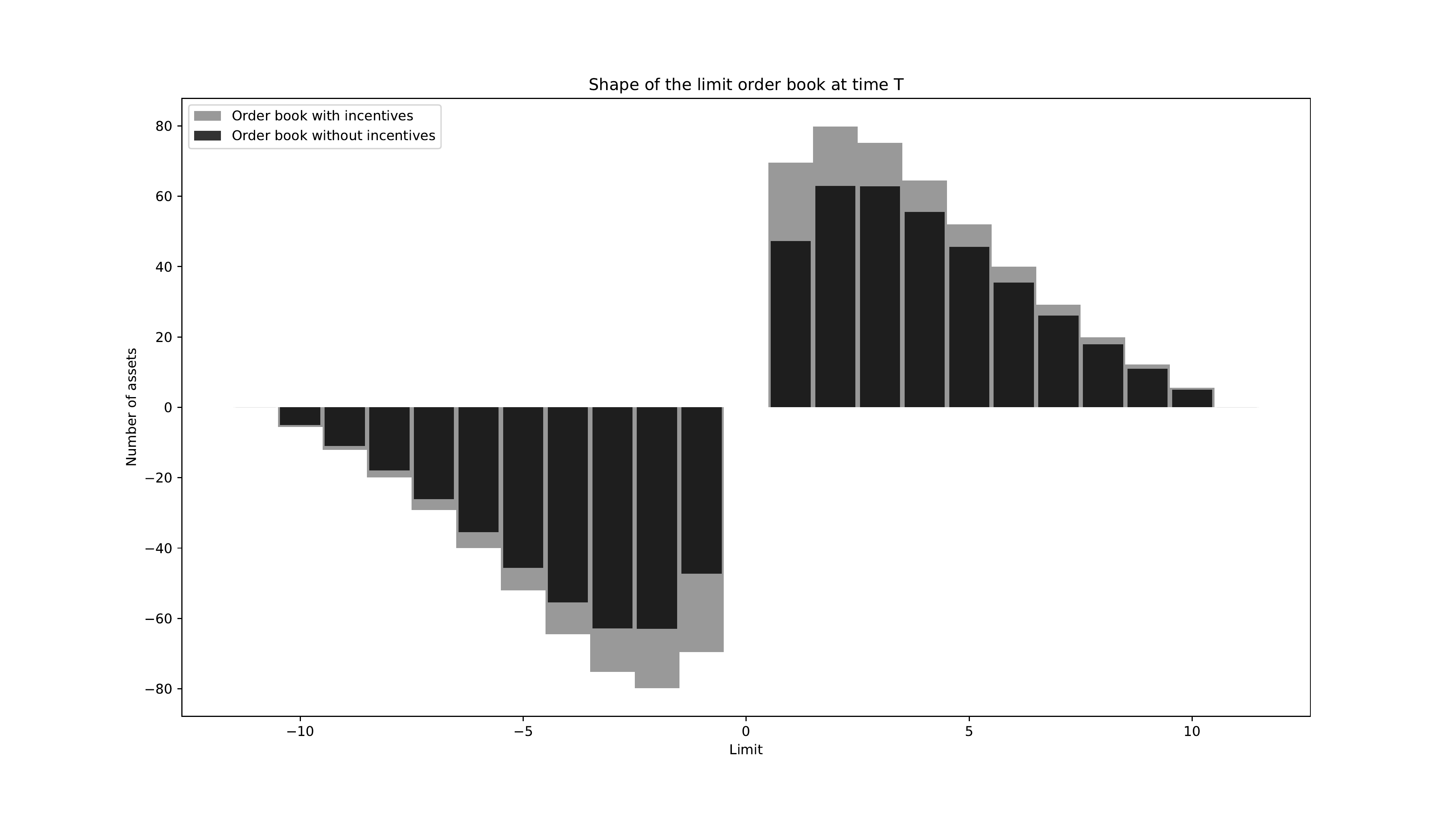}\\
\caption{Shape of the limit order book at time $T$.}\label{shape_3}
\end{figure}

We observe that an increase in the convection rate also leads to an increase in the incentives. Indeed, when $\alpha_b$ and $\alpha_a$ increase, the exchange needs to compensate the higher cancellation rate by setting more attractive incentives. This allows to get roughly the same level and shape for the order book with incentives in the case with high cancellation rate plotted in Figure \ref{shape_3} as in the first case in Figure \ref{shape_0}.

\newpage 

\section{Conclusion}
In this paper, we tackled the problem of incentives design in a limit order book. We represented this problem as an optimal control problem for a SPDE, using the limit order book model of Cont and Müller \cite{cont2021stochastic}. Due to the specific form of the control problem faced by the exchange, we showed that the solution can be expressed with a Feynman-Kac representation. Moreover, we gave an example of a class of intensity and penalty functions under which the asymptotic behavior of the control problem leads to closed-form incentives at each limit of the order book. We then simulated a fictitious limit order book using these incentives to illustrate numerically their influence on its shape. In particular, we showed that these incentives increase liquidity at the first limits of the order book.  Finally, we performed a sensitivity analysis with respect to the market parameters. 
\appendix 

\section{Appendix}\label{sec_proofs}

\subsection{Proof of Theorem \ref{maximum}}\label{Proof1}

We follow closely the proof of Theorem 10.2 in \cite{oksendal2007applied}.\\

Let $Z^a \in \mathcal A^a$, and for the sake of simplicity, denote by $u$ the associated solution to SPDE \eqref{SPDE1} : $u:= u^{Z^a}$. For $(t,x) \in [0,T) \times (0,L)$, we introduce
\begin{align*}
\hat H_t(x) &:= H\left(x, \hat u_t(x), \hat Z^a_t(x), \hat p_t(x), \hat q_t(x) \right)\\
&=\hat u_t(x) - g(x,\hat Z^a_t(x)) + \left(\alpha_a \hat u_t(x) +f^a(x,\hat Z^a_t(x))\right)\hat p_t(x) + \sigma_a \hat u_t(x) \hat q_t(x),
\end{align*}
and 
\begin{align*}
H_t(x) &:= H\left(x, u_t(x), Z^a_t(x), \hat p_t(x), \hat q_t(x) \right)\\
&= u_t(x) - g(x, Z^a_t(x)) + \left(\vphantom{\hat Z^a_t(x)}\alpha_a  u_t(x) +f^a(x,Z^a_t(x))\right)\hat p_t(x) + \sigma_a  u_t(x) \hat q_t(x).
\end{align*}
Hence we have 
$$J(\hat Z^a) - J(Z^a) = I,$$
where 
\begin{align*}
I&:= \mathbb E \left[ \int_0^T \int_0^L \left \{ \hat u_t(x) - u_t(x) -\left(g(x,\hat Z^a_t(x)) - g(x,Z^a_t(x)) \right)  \right\} dx dt \right]\\
&= \mathbb E\! \bigg[ \int_0^T\!\!\! \int_0^L\!\! \Big\{ \hat H_t(x)\! -\! H_t(x)\! -\!\big(\alpha_a (\hat u_t(x) - u_t(x)) + f^a(x,\hat Z^a_t(x))- f^a(x,Z^a_t(x)) \big) \hat p_t(x) \\
& \qquad \qquad \qquad \qquad \qquad \qquad \qquad \qquad \qquad \qquad \qquad \qquad \qquad \qquad   - \sigma_a (\hat u_t(x)\! -\! u_t(x))\hat q_t(x)  \Big\} dx dt \bigg].
\end{align*}
Moreover, denoting by $\tilde u_t(x)$ the process 
$$\tilde u_t(x) := u_t(x) - \hat u_t(x),$$
we know from the boundary condition \eqref{boundBSPDE} that 
\begin{align*}
0 &= \mathbb E\left[\int_0^L \hat p_T(x)\tilde u_T(x) dx   \right]\\
&= \mathbb E\left[\int_0^L \left( \hat p_0(x)\tilde u_0(x)  + \int_0^T \left\{ \tilde u_t(x) d\hat p_t(x) + \hat p_t(x) d\tilde u_t(x) +  \sigma_a \tilde u_t(x) \hat q_t(x) dt\right\}   \right) dx\right]\\
&= \mathbb E \Bigg[ \int_0^L \int_0^T \bigg\{ \tilde u_t(x) \left(- \mathcal L^* \hat p_t(x) - \widehat{\partial_u H}_t(x) \right)\\
& \qquad + \hat p_t(x) \left( \mathcal L \tilde u_t(x) +\left(\alpha_a \tilde u_t(x) +f^a(x,Z^a_t(x)) - f^a(x, \hat Z^a_t(x))\right) \right) + \sigma_a \tilde u_t(x) \hat q_t(x) \bigg\}dt dx \Bigg],
\end{align*}
where
$$\widehat{\partial_u H}_t(x) = \partial_u H \left(x, \hat u_t(x) \hat Z^a_t(x), \hat p_t(x), \hat q_t(x) \right) \quad \forall (t,x) \in [0,T) \times (0,L).$$
Hence 
\begin{align*}
J(\hat Z^a) - J(Z^a)& = I - \mathbb E\left[\int_0^L \hat p_T(x)\tilde u_T(x) dx   \right] \\
&= \mathbb E \Bigg[ \int_0^T\left( \int_0^L \left\{ \tilde u_t(x) \mathcal L^* \hat p_t(x) - \hat p_t(x) \mathcal L \tilde u_t(x) \right\} dx\right)dt\\
& \qquad + \int_0^L \left( \int_0^T \left\{ \hat H_t(x) - H_t(x) + \widehat{\partial_u H}_t(x) \tilde u_t(x) \right\} dt \right) dx  \Bigg].
\end{align*}
As $\tilde u_t(x) = \hat p_t(x) = 0 $ $ \forall (t,x) \in [0,T) \times \{0,L\}$ a.s, we have a.s. 
$$\int_0^L  \tilde u_t(x) \mathcal L^* \hat p_t(x)    dx = \int_0^L \hat p_t(x) \mathcal L \tilde u_t(x) dx, $$
and therefore
$$J(\hat Z^a) - J(Z^a) = \mathbb E \Bigg[ \int_0^L \left( \int_0^T \left\{ \hat H_t(x) - H_t(x) + \widehat{\partial_u H}_t(x) \tilde u_t(x) \right\} dt \right) dx  \Bigg].$$
Finally by Lemma \ref{Hconcave}, we have
$$H_t(x) - \hat H_t(x) \le \widehat{\partial_u H}_t(x) \tilde u_t(x) + \widehat{\partial_z H}_t(x) \left(Z^a_t(x) - \hat Z^a_t(x) \right),$$
where
$$\widehat{\partial_z H}_t(x) = \partial_z H \left(x, \hat u_t(x) \hat Z^a_t(x), \hat p_t(x), \hat q_t(x) \right) \quad \forall (t,x) \in [0,T) \times (0,L).$$
But by $(i),$ 
$$\widehat{\partial_z H}_t(x) \left(Z^a_t(x) - \hat Z^a_t(x) \right) \le 0,$$
and therefore
$$\hat H_t(x) - H_t(x) + \widehat{\partial_u H}_t(x) \tilde u_t(x) \ge 0,$$
and
$$J(\hat Z^a) \ge J(Z^a).$$

\bibliographystyle{plainnat}
\bibliography{biblio.bib}

\end{document}